\documentclass{amsart}
\usepackage{amsmath}
\usepackage{amssymb}
\usepackage{amscd}
\usepackage{latexsym}
\usepackage{graphicx} 

\theoremstyle{plain}
\newtheorem{thm}{Theorem}[section]

\newtheorem{prop}[thm]{Proposition}

\newtheorem*{ratthm}{Rational Blowdown Theorem}
\newtheorem*{ksthm}{Knot Surgery Theorem}

\theoremstyle{definition}

\def \CPb {\overline{\mathbf{CP}}^{2}}
\def \CP {{\mathbf{CP}}^{2}} 
\def \R {\mathbf{R}}
\def \Z {\mathbf{Z}}
\def \Sig{\Sigma}

\def \vp {\varphi}

\def \lam {\lambda}
\def\L{\Lambda}
\def \G {\Gamma}
\def \l {\ell}
\def \o {\omega}

\def \bd {\partial}
\def \x {\times}
\def \- {\!\smallsetminus\!}
\def \C {\subset}

\def \ssw {\text{SW}}
\def \sw {\mathcal{SW}}

\def \DD {\Delta}

\def \Ds {\DD^{\! s}}

\def\hk{\widehat{k}}

\begin{document}

\baselineskip.5cm
\title [Surgery on nullhomologous tori] {Surgery on nullhomologous tori and simply connected 4-manifolds with ${\bf{b^+=1}}$}
\author[Ronald Fintushel]{Ronald Fintushel}
\address{Department of Mathematics, Michigan State University \newline
\hspace*{.375in}East Lansing, Michigan 48824}
\email{\rm{ronfint@math.msu.edu}}
\thanks{The first author was partially supported NSF Grant DMS0305818
and the second author by NSF Grant DMS0505080}
\author[Ronald J. Stern]{Ronald J. Stern}
\address{Department of Mathematics, University of California \newline
\hspace*{.375in}Irvine,  California 92697}
\email{\rm{rstern@math.uci.edu}}

\maketitle

\section{Introduction\label{Intro}}

In the past few years there has been significant progress on the problem of finding exotic smooth structures on the manifolds $P_m=\CP\# m\,\CPb$. The initial step was taken by Jongil Park, \cite{P}, who found the first exotic smooth structure on $P_7$, and whose ideas renewed the interest in this subject. Peter Ozsvath and Zoltan Szabo proved that Park's manifold is minimal \cite{OS}, and Andras Stipsicz and Szabo used a technique similar to Park's to construct an exotic structure on $P_6$ \cite{SS}. Shortly thereafter, the authors of this article developed a technique for producing infinite families of smooth structures on $P_m$, $6\le m\le 8$ \cite{DN}, and Park, Stipsicz, and Szabo showed that this can be applied to the case $m=5$ \cite{PSS}.

It is the goal of this paper to better understand the underlying mechanism which produces infinitely many distinct smooth structures on $P_m$, $5 \le m\le 8$. As we explain below, all these constructions start with the elliptic surface $E(1)=P_9$, perform a knot surgery using a family of twist knots indexed by an integer $n$ \cite{KL4M},  then blow the result up several times in order to find a suitable configuration of spheres that can be rationally blown down \cite{rat} to obtain a smooth structure on $P_{m}$ that is distinguished by the integer $n$. We shall explain how this can be accomplished by surgery on nullhomologous tori in a manifold $R_m$ homeomorphic to $P_m$, $5\le m\le 8$. In other words, we shall find a nullhomologous torus $\L_m$ in $R_m$ so that $1/n$-surgery on $\L_m$ preserves the homeomorphism type of $R_m$, but changes the smooth structure of $R_m$ in a way that depends on $n$. Presumably, $R_m$ is diffeomorphic to $P_m$, but we have not yet been able to show this in general. Our hope is that by better understanding $\L_m$ and its properties, one will be able to find similar nullhomologous tori in $P_m$, for $m<5$. 

\section{A short history of simply connected $4$-manifolds with $b^+=1$}

It is a basic question of $4$-manifold topology to understand the smooth structures on the complex projective plane $\CP$. Thus one is interested in knowing the smallest 
$m$ for which $P_m=\CP\# m\,\CPb$ admits an exotic smooth structure. The first such example was produced by Simon Donaldson in the historic paper \cite{D}, where it was shown that the Dolgachev surface $E(1)_{2,3}$, the result of performing log transforms of orders 2 and 3 on the rational elliptic surface $E(1)=P_9$ \cite{Dv}, is homeomorphic but not diffeomorphic to $P_9$. This breakthrough example provided the first known instance of an exotic smooth structure on a simply connected $4$-manifold. Later, work of Friedman and Morgan \cite{FM} showed that the integers $p$ and $q$ are smooth invariants. The fact that $E(1)_{p,q}$ is not diffeomorphic to $E(1)_{p',q'}$ for $\{ p,q\}\ne\{ p',q'\}$ persists even after an arbitrary number of blowups \cite{blowup}; however, no minimal exotic smooth structures are currently known for $P_m$, $m\ge 10$. 

In the late 1980's, Dieter Kotschick \cite{K} proved that the Barlow surface, known to be homeomorphic to $P_8$, is not diffeomorphic to it. However,
in following years the subject of simply connected smooth $4$-manifolds with $b^+=1$  languished because of a lack of suitable examples.  As we mentioned above, largely due to the example of Park \cite{P}  of an exotic smooth structure on $P_7$, this topic has again become active. 

Here is an outline of a version of Park's example: Consider $E(1)$ with an elliptic fibration whose singular fibers are four nodal fibers and and $I_8$-fiber. (An $I_n$-fiber is comprised of a circular plumbing of $n$ $2$-spheres of self-intersection $-2$ \cite{BPV}.) This elliptic fibration has a section which is an exceptional curve $E$ (of self-intersection $-1$). Blow up $E(1)$ four times, at the double points of the four nodal fibers. Then in 
$E(1)\#4\CPb \cong P_{13}$, we find a configuration of $2$-spheres consisting of $E$, four disjoint spheres of self-intersection $-4$, each intersecting $E$ once, and the $I_8$-fiber, which $E$ intersects in exactly one $2$-sphere, and the $I_8$-fiber is disjoint from the four spheres of self-intersection $-4$.

The transverse intersections of $E$ with the four spheres of self-intersection $-4$ can be smoothed to obtain a 2-sphere of self-intersection $-9$. Together with spheres from the $I_8$-fiber, we obtain a linear configuration of $2$-spheres:

\vspace{-.2in}\hspace{.25in}
\begin{picture}(100,60)(-90,-25)
 \put(-5,3){\makebox(200,20)[bl]{$\ -9$ \hspace{3pt}
                                  $-2$\hspace{7pt}$-2$\hspace{6pt}$-2$\hspace{7pt}$-2$\hspace{7pt}$-2$}}
  \multiput(10,0)(20,0){5}{\line(1,0){20}}
  \multiput(10,0)(20,0){6}{\circle*{3}}
\end{picture}

\vspace{-.2in}\noindent A regular neighborhood of this configuration has as its boundary the lens space $L(49,-6)$, and as we explain below, this lens space bounds a rational homology ball. This means that this configuration can be rationally blown down \cite{rat}, reducing $b^-$ by $6$. One obtains a manifold $P$ with $b^+=1$ and $b^-=7$. It is not difficult to show that $P$ is simply connected, and so it is homeomorphic to $P_7$. It follows from Seiberg-Witten theory that $P$ is not diffeomorphic to $P_7$. This will be explained below. While this is not precisely the description of the manifold that was given in \cite{P} (and it is not even clear that $P$ is diffeomorphic to the example of \cite{P}), the construction given here is similar to that of Park.

Stipsicz and Szabo improved on Park's example by finding a more complicated configuration in a larger blowup of $E(1)$, yet one which could be rationally blown down to get an even smaller manifold, homeomorphic but not diffeomorphic to $P_6$. 

For some time, even after these examples, it was suspected that $P_m$, $m < 9$ would support only finitely many distinct smooth structures. This was due to the fact that until \cite{DN} the only technique available for producing infinitely many distinct smooth structures on a given smooth $4$-manifold $X$ was to require that $X$ contain a minimal genus torus with trivial normal bundle and representing a nontrivial homology class. It is known that $P_m$, $m < 9$, contain no such tori. Thus it is is the goal of this paper to better understand techniques for producing infinitely many distinct smooth structures and to better understand the examples of  \cite{DN}.

\section{Seiberg-Witten invariants, rational blowdowns, and knot surgery}

\subsection{Seiberg-Witten invariants} Let $X$ be a simply connected oriented 
4-mani-fold with $b_X^+=1$ with a given orientation of $H^2_+(X;\R)$ and a given metric $g$. The Seiberg-Witten invariant depends on the metric $g$ and a self-dual 2-form as follows. 
There is a unique $g$-self-dual harmonic 2-form $\o_g\in H^2_+(X;\R)$ with $\o_g^2=1$ and corresponding to the positive orientation. (Often $\o_g$ is called a {\it{period point}} for the metric $g$.) Fix a characteristic homology class $k\in H_2(X;\Z)$.  Given a pair $(A,\psi)$, where
$A$ is a connection in the complex line bundle whose first Chern class is the Poincar\'e dual $\hk=\frac{i}{2\pi}[F_A]$ of $k$ and $\psi$ a section of the bundle $W^+$ of self-dual spinors for the associated $spin^{\,c}$ structure, the perturbed Seiberg-Witten equations are:
\begin{gather*} 
D_A\psi = 0 \\
F_A^+  = q(\psi)+i\eta \notag\label{SWeqn}
\end{gather*}
where $F_A^+$ is the self-dual part of the curvature $F_A$,
$D_A$ is the twisted Dirac operator, $\eta$ is a
self-dual 2-form on $X$, and
$q$ is a quadratic function. Write $\ssw_{X,g,\eta}(k)$ for the
corresponding invariant. As the pair
$(g,\eta)$ varies, $\ssw_{X,g,\eta}(k)$ can change only at those pairs
$(g,\eta)$ for which there are solutions with $\psi=0$. These 
solutions occur for pairs $(g,\eta)$ satisfying $(2\pi\hk+\eta)\cdot\o_g=0$.
This last equation defines a wall in $H^2(X;\R)$. 

The point $\o_g$ determines a component of the double cone consisting of elements of $H^2(X;\R)$ of positive square. We prefer to work with $H_2(X;\R)$. The dual component is determined by the Poincar\'e dual $H$ of $\omega_g$.  An element $H'\in H_2(X;\R)$ of positive square lies in the same component as $H$ if $H'\cdot H>0$. If
$(2\pi \hk+\eta)\cdot\o_g\ne 0$ for a generic $\eta$, $\,\ssw_{X,g,\eta}(k)$ is
well-defined, and its value depends only on the sign of $(2\pi \hk+\eta)\cdot\o_g$. Write $\ssw_{X,H}^+(k)$ for $\ssw_{X,g,\eta}(k)$ if 
$(2\pi \hk+\eta)\cdot\o_g>0$ and $\ssw_{X,H}^-(k)$ in the other case.

The invariant $\ssw_{X,H}(k)$ is defined by $\ssw_{X,H}(k) =\ssw_{X,H}^+(k)$ if 
$(2\pi \hk)\cdot\o_g>0$, or dually, if $k\cdot H>0$, and $\ssw_{X,H}(k) =\ssw_{X,H}^-(k)$ if $k\cdot H <0$. The wall-crossing formula \cite{KM,LL} states that if $H', H''$ are elements of positive square in $H_2(X;\R)$ with $H'\cdot H>0$ and $H''\cdot H>0$, then if $k\cdot H' <0$ and $k\cdot H''>0$,
\[ \ssw_{X,H''}(k) - \ssw_{X,H'}(k) = (-1)^{1+\frac12 d(k)}\]
where $d(k)=\frac14(k^2-(3\,\text{sign}+2\,\text{e})(X))$ is the formal dimension of the Seiberg-Witten moduli spaces.

Furthermore, in case $b^-\le 9$, the wall-crossing formula, together with the fact that $\ssw_{X,H}(k)=0$ if $d(k)<0$, implies that $\ssw_{X,H}(k) = \ssw_{X,H'}(k)$ for any $H'$ of positive square in $H_2(X;\R)$ with $H\cdot H'>0$. So in case $b^+_X=1$ and $b^-_X\le 9$, there is a well-defined Seiberg-Witten invariant, $\ssw_X(k)$. If $\ssw_X(k)\ne 0$, $k$ is called a {\it{basic class}} of $X$.

It is convenient to view the Seiberg-Witten invariant as an element of the integral group ring $\Z H_2(X)$. For $k\in H_2(X)$ we let $t_k$ denote the corresponding element in $\Z H_2(X)$. Then the Seiberg-Witten invariant of $X$ is
\[\sw_{X,H} = \sum\ssw_{X,H}(k)\cdot t_k\]

An important property of the Seiberg-Witten invariant is that if $X$ admits a metric $g$ of positive scalar curvature, then for the Poincar\'e dual $H$ of $\omega_g$, we have 
$\sw_{X,H}=0$. In particular, for $m\le 9$, $\sw_{P_m}=0$.

\subsection{Rational blowdowns} Let $C_{p}$ be the smooth $4$-manifold obtained by plumbing $(p-1)$ disk bundles over the $2$-sphere according to
the diagram

\begin{picture}(100,60)(-90,-25)
 \put(-12,3){\makebox(200,20)[bl]{$-(p+2)$ \hspace{6pt}
                                  $-2$ \hspace{96pt} $-2$}}
 \put(4,-25){\makebox(200,20)[tl]{$u_{0}$ \hspace{25pt}
                                  $u_{1}$ \hspace{86pt} $u_{p-2}$}}
  \multiput(10,0)(40,0){2}{\line(1,0){40}}
  \multiput(10,0)(40,0){2}{\circle*{3}}
  \multiput(100,0)(5,0){4}{\makebox(0,0){$\cdots$}}
  \put(125,0){\line(1,0){40}}
  \put(165,0){\circle*{3}}
\end{picture}

\noindent Then the classes of the $0$-sections have self-intersections $u_0^2=-(p+2)$ and $u_i^2=-2$, $i=1,\dots,p-2$. The boundary of $C_p$ is the lens space 
$L(p^2, 1-p)$ which bounds a rational ball $B_p$ with $\pi_1(B_p)=\Z_p$ and $\pi_1(\bd B_p)\to \pi_1(B_p)$ surjective. If $C_p$ is embedded in a $4$-manifold $X$ then the rational blowdown manifold $X_{(p)}$ of  \cite{rat} is obtained by replacing $C_p$ with $B_p$, i.e., $X_{(p)} = (X\- C_p) \cup B_p$. This construction is independent of the choice of gluing map.

\begin{ratthm}[\cite{rat}] Let $X$ be a simply connected $4$-manifold containing the configuration of $2$-spheres, $C_p$. If $X\- C_p$ is also simply connected then so is the rational blowdown manifold $X_{(p)}$. The homology, $H_2(X_{(p)};\R)$, may be identified with the orthogonal complement of the classes $u_i$, $i=0,\dots,p-2$ in $H_2(X;\R)$. 

Given a characteristic homology class $k\in H_2(X_{(p)};\Z)$, there is a lift $\widetilde{k} \in H_2(X;\Z)$ which is characteristic and for which the dimensions of moduli spaces agree, $d_{X_{(p)}}(k)=d_X(\widetilde{k})$. If $b^+_X>1$ then $\ssw_{X_{(p)}}(k)=\ssw_X(\widetilde{k})$. In case $b^+_X=1$, if $H\in H_2^+(X;\R)$ is orthogonal to all the $u_i$ then it also can be viewed as an element of $H_2^+(X_{(p)};\R)$, and 
$\ssw_{X_{(p)},H}(k)=\ssw_{X,H}(\widetilde{k})$.
\end{ratthm}

\subsection{Knot surgery} Let $X$ be a $4$-manifold which contains a homologically essential torus $T$ of self-intersection $0$, and let $K$ be a knot in $S^3$. Let $N(K)$ be a tubular neighborhood of $K$ in $S^3$, and let $T\x D^2$ be a tubular neighborhood of $T$ in $X$. Then the knot surgery manifold $X_K$ is defined by 
\[ X_K = (X\- (T\x D^2))\cup (S^1\x (S^3\- N(K))\]
The two pieces are glued together in such a way that 
the homology class $[{\text{pt}}\x \bd D^2]$ is identified with $[{\text{pt}}\x \lam]$ where $\lam$ is the class of a longitude of $K$. For example, if X is a simply connected elliptic surface with a (spherical) section $S$ of self-intersection $n$ and one performs knot surgery on the fiber $T$ of this fibration, then the gluing condition implies that in $X_K$ there is a pseudosection $S_K$ of genus equal to the genus of the knot $K$ and with self-intersection $n$. By `pseudosection' we mean that the intersection number $S_K\cdot F=1$. However, $X_K$  need no longer be an elliptic surface. This surface $S_K$ is constructed by removing a disk from $S$ where it intersects the fiber $T$ and replacing this disk by a Seifert surface for the knot $K$.

One can also interpret $X_K$ as a fiber sum. Let $M_K$ denote the $3$-manifold obtained from $0$-framed surgery on $K$ in $S^3$. Then 
\[ X_K = X \#_{T=S^1\x m} S^1\x M_K \]
where $m$ is a meridian of $K$.

The gluing
condition does not, in general, completely determine the diffeomorphism type of $X_K$; however if we take $X_K$ to be any manifold constructed in this fashion and if, for example, $T$ has a cusp neighborhood, then the Seiberg-Witten invariant of $X_K$ is completely determined by the Seiberg-Witten invariant of $X$ and the symmetrized Alexander polynomial $\Ds$ of $K$:

\begin{ksthm}[\cite{KL4M}] Let $X$ be a $4$-manifold which contains a homologically essential torus $T$ of self-intersection $0$ whose $H_1$ is generated by vanishing cycles, and let $K$ be a knot in $S^3$. The Seiberg-Witten invariant of the knot surgery manifold $X_K$ is given by
 \[ \sw_{X_K}=\sw_X\cdot\Ds(t^2) \]
where $t$ represents the homology class of the torus $T$. Furthermore, if $X$ and $X\- T$ are simply connected, then so is $X_K$. 
\end{ksthm}

\section{Double node neighborhoods and knot surgery}

A simply connected elliptic surface is fibered over $S^2$ with smooth fiber a torus and with singular fibers. The most generic type of singular fiber is a nodal fiber (an immersed $2$-sphere with one transverse positive double point). The monodromy of a nodal fiber is $D_a$, a Dehn twist around the `vanishing cycle' $a\in H_1(F;\Z)$, where $F$ is a smooth fiber of the elliptic fibration. The vanishing cycle $a$ is represented by a nonseparating loop on the smooth fiber and the nodal fiber is obtained by collapsing this vanishing cycle to a point to create a transverse self-intersection. The vanishing cycle  bounds a `vanishing disk', a disk of relative self-intersection $-1$ with respect to the framing of its boundary given by pushing the loop off itself on the smooth fiber. 

An $I_2$-fiber consists of a pair of $2$-spheres of self-intersection $-2$ which are plumbed at two points. The monodromy of an $I_2$-fiber is $D_a^2$, which is also the monodromy of a pair of nodal fibers with the same vanishing cycle. This means that an elliptic fibration which contains an $I_2$-fiber can be perturbed to contain two nodal fibers with the same vanishing cycle.  

A {\em{double node neighborhood}} $D$ is a fibered neighborhood of an elliptic fibration which contains exactly two nodal fibers with the same vanishing cycle. If $F$ is a smooth fiber of $D$, there is a vanishing class $a$ that bounds vanishing disks in the two different nodal fibers, and these give rise to a sphere $V$ of self-intersection $-2$ in $D$. 

In \cite{DN} we showed how performing knot surgery in a double node neighborhood $D$ in $E(1)$ can give rise to an immersed pseudosection of self-intersection $-1$. Let us review a version of this construction. Consider the knot $K$ of Figure 1. Of course, this is just the unknot, and we see a Seifert surface $\Sig$ of genus one. Let $\G$ be the loop which runs through both half-twists in the clasp. Then $\G$ satisfies the two key conditions of \cite{DN}:
\begin{itemize}
\item[(i)] $\G$ bounds a disk in $S^3$ which intersects $K$ in exactly two points.
\item[(ii)] The linking number in $S^3$ of $\G$ with its pushoff on $\Sig$ is $+1$.
\end{itemize}   
It follows from these properties that $\G$ bounds a punctured torus in $S^3\- K$. 

\smallskip

\hspace{1.6in}\includegraphics[scale=0.50]{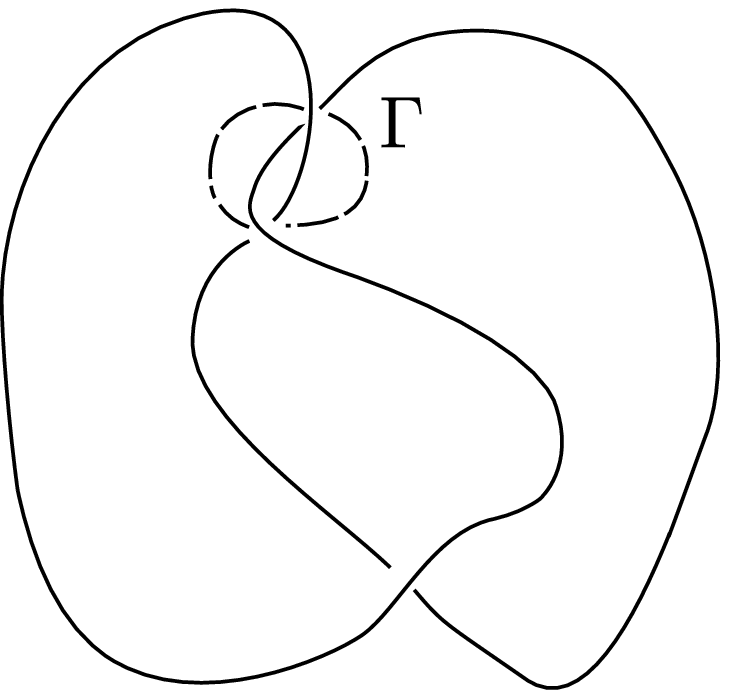}\vspace{-.2in}\vspace{.25in}\\
\centerline { Figure 1: $K =$ unknot}

\smallskip
It is known that $E(1)$ admits an elliptic fibration with two nodal fibers, an $I_2$-fiber and an $I_8$-fiber \cite{Pn}. As above, this fibration can be perturbed so that the $I_2$-fiber gives us a double node neighborhood $D$ with vanishing cycle $a$. 
Consider the result of knot surgery in $D$ using the knot $K$ and the fiber $F$ of $E(1)$. In the knot surgery construction, one is free to make any choice of gluing as long as a longitude of $K$ is sent to the boundary circle of a normal disk to $F$. We choose the gluing so that the class of a meridian $m$ of $K$ is sent to the class of $a \x \{pt\}$ in $H_1(\bd(D\- N(F));\Z)=H_1(F\x \bd D^2;\Z)$. Note that the result of knot surgery, 
\[ E(1)_K= E(1)\- N(F) \cup S^1\x(S^3\- N(K)) =  E(1)\- N(F) \cup T^2\x D^2 \]
because $K$ is the unknot. Since any diffeomorphism of $\bd (E(1)\- N(F))$ extends over all of $E(1)\- N(F)$, we see that $E(1)_K$ is diffeomorphic to $E(1)$.

There is a genus one pseudosection $S_K$ in $E(1)_K$ which is formed using the genus one Seifert surface $\Sig$. The self-intersection of $S_K$ is $-1$.  The loop $\G$ sits on $S_K$ and by (i) it bounds a twice-punctured disk $\DD$ in $\{pt\}\x \bd(S^3\- N(K))$ where $\bd\DD=\G\cup m_1\cup m_2$ where the $m_i$ are meridians of $K$. The meridians $m_i$ bound disjoint vanishing disks $\DD_i$ in $D\- N(F)$ since they are identified with disjoint loops each of which represents the class of $a \x \{pt\}$ in $H_1(\bd(D\- N(F));\Z)$. Hence  in $D_K$, the result of knot surgery on $D$,  the loop $\G\C S_K$ bounds a disk $U=\DD\cup \DD_1\cup \DD_2$. By construction, the relative self-intersection of $U$ relative to the framing given by the pushoff of $\G$ in $S_K$ is $+1-1-1=-1$. (This uses (ii).) Furthermore, $U\cap S_K=\G$.

Since $\G$ is nonseparating in $S_K$, surgery on it kills $\pi_1(S_K)$. Ambient surgery may be performed in $D_K$ by removing an annular neighborhood of $\G$ and replacing it with a pair of disks $U'$, $U''$ as obtained above. These disks intersect in a single point, and this is precisely the complex-algebraic model of a nodal intersection. This means that we can represent the homology class of the pseudosection $[S_K]$ in $H_2(E(1)_K;\Z)$ by an immersed sphere $S'$ with one positive double point.

With these as preliminaries, our goal for the remainder of this paper is to construct for for every $5 \le m \le 8$ a manifold $R_{m}$ that is homeomorphic to $P_{m}$, has vanishing Seiberg-Witten invariants, and contains a nullhomologous torus $T_{m}$ with the property that a $1/n$-surgery on $T_{m}$ (with respect to the nullhomologous framing) yields a smooth structure on $P_{m}$ distinguished by the integer $n$. We conjecture that $R_{m}$ is diffeomorphic to $P_{m}$, but we are unable to show this at this time. We start with the $m=8$ case in the next section.

\section{Infinite families homeomorphic to $\CP\#8\CPb$}

The construction of the previous section shows that in $E(1)\cong E(1)_K$ one has the configuration consisting of the immersed $2$-sphere $S'$ with a pair of disjoint nodal fibers, each intersecting $S'$ once transversely. Also, $S'$ intersects the $I_8$-fiber transversely in one point. This is illustrated in Figure 2.

At this stage there are three possibilities:

\noindent 1.\ If we blow up the double point of $S'$, then in $P_{10}$ we obtain a configuration consisting of the total transform $S''$ of $S'$, which is a sphere of self-intersection $-5$, and the sphere of self-intersection $-2$ at which $S'$ intersects $I_8$.  (see Figure 3.) This is the configuration $C_3$ which can be rationally blown down to obtain a manifold $R_8$ with $b^+=1$ and $b^-=8$. It is easy to see that $R_8$ is simply connected; so $R_8$ is homeomorphic to $P_8$.

\bigskip

\hspace{.5in}\includegraphics[scale=0.50]{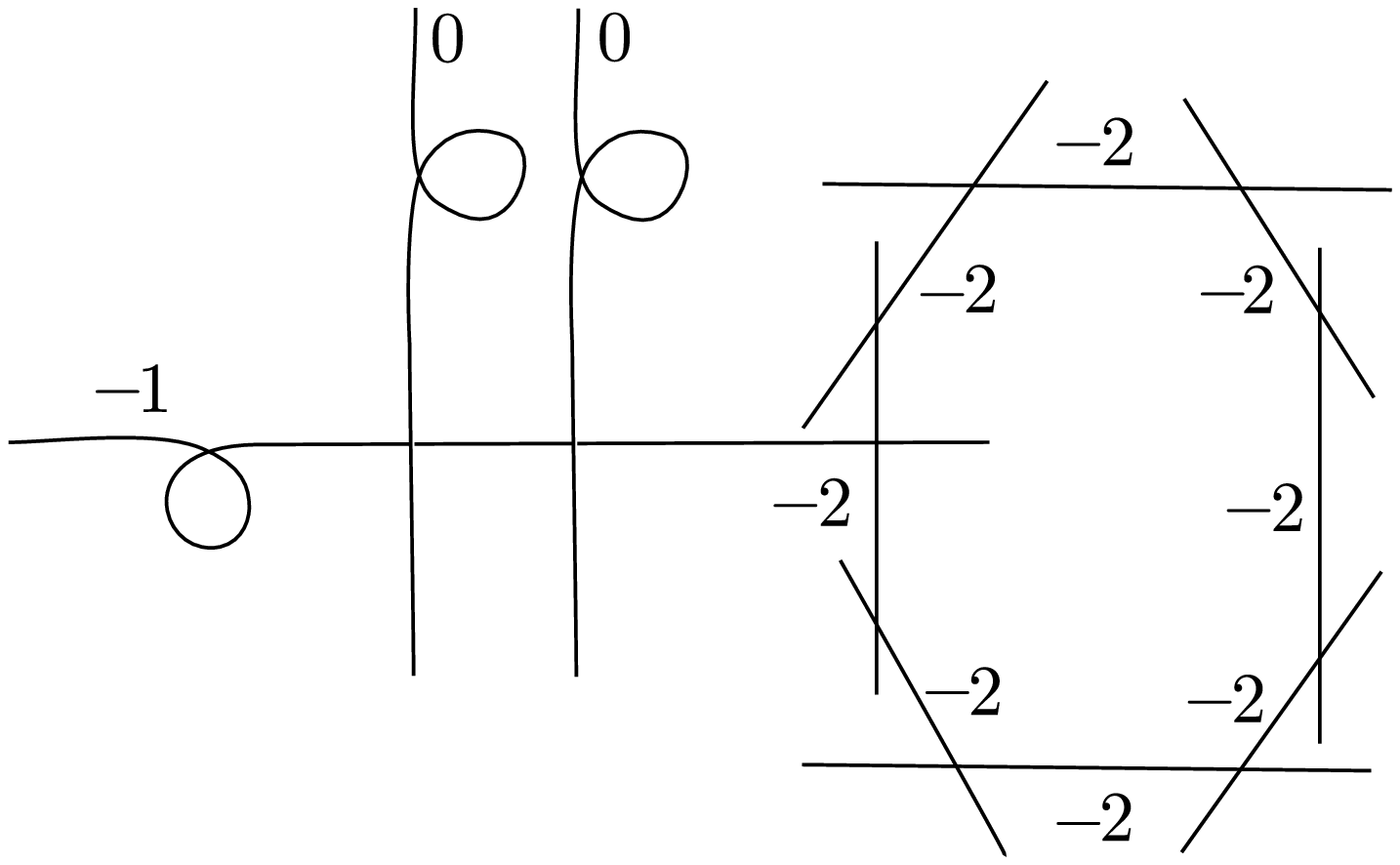}\\
\centerline{{Figure 2}}

\bigskip

\hspace{.5in}\includegraphics[scale=0.50]{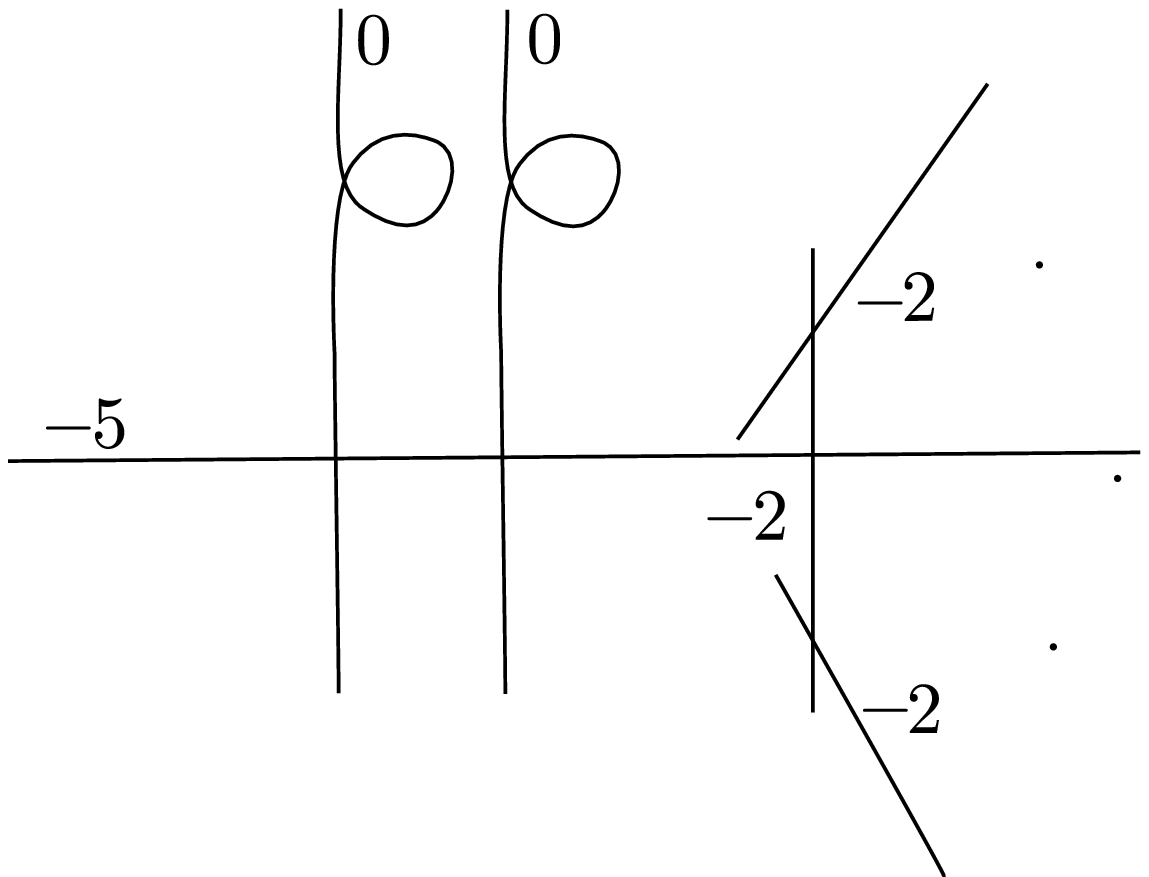}\\
\centerline{{Figure 3}}

\medskip

\noindent 2.\ Blow up at the double point of $S'$ as well as at the double point of one of the nodal fibers. Then in $P_{11}$ we get a configuration of $2$-spheres consisting of $S''$, a transverse sphere $F'$ of self-intersection $-4$, and three spheres from the $I_8$-fiber.  (See Figure~4.) Smoothing the intersection of  $S''$ and $F'$ gives a sphere of self-intersection $-7$ and we obtain the configuration $C_5$ in $P_{11}$. Rationally blowing down $C_5$ gives a manifold $R_7$ homeomorphic to $P_7$.

\hspace{.4in}\includegraphics[scale=0.50]{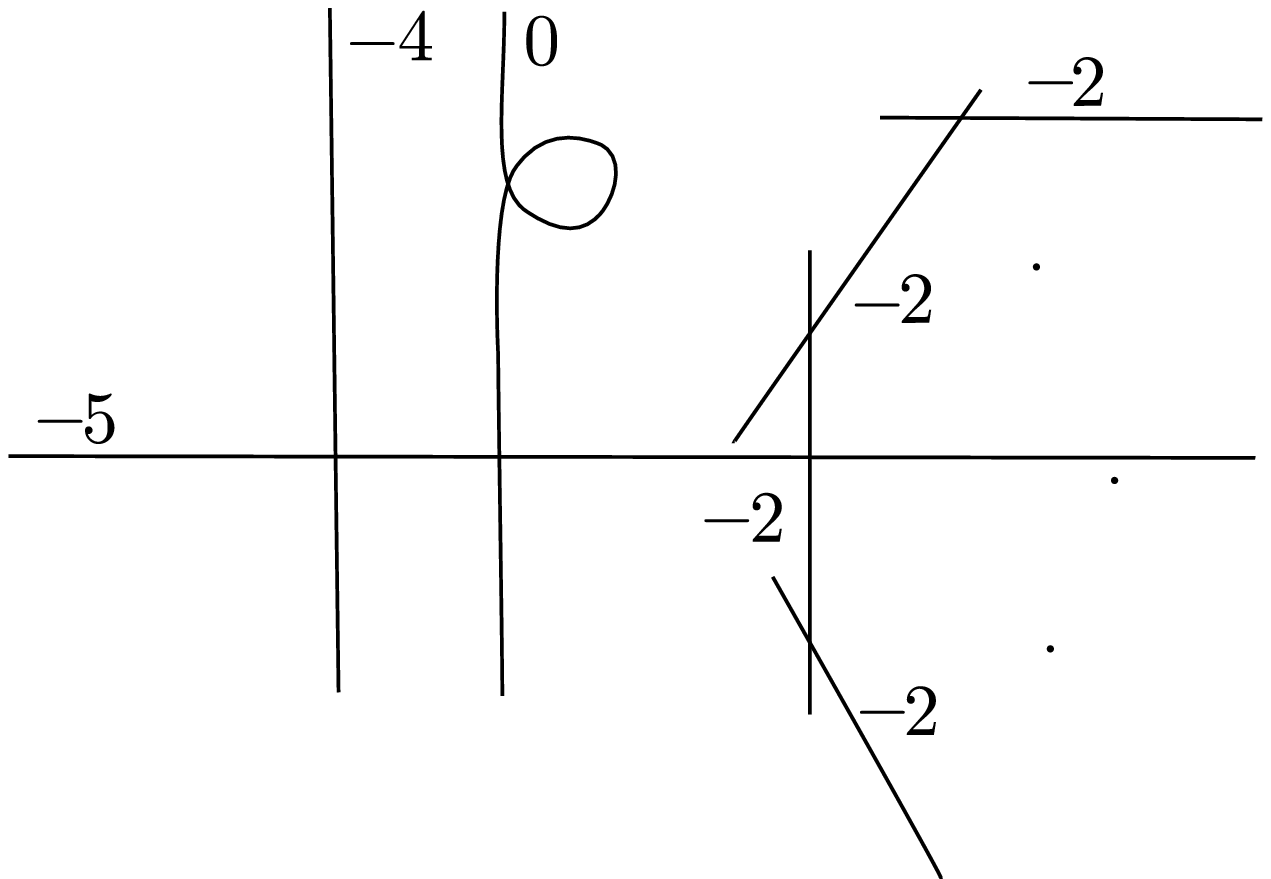}\\
\centerline{{Figure 4}}

\medskip

\vspace{.1in}\noindent 3.\  Blow up at the double point of $S'$ as well as at the double points of both nodal fibers. Then in $P_{12}$ we get a configuration of $2$-spheres consisting of $S''$, two disjoint transverse spheres $F'$, $F''$ of self-intersection $-4$, and five spheres from the $I_8$-fiber. (See Figure~5.) Smoothing the intersections of  $S''$, $F'$ and $F''$ gives a sphere of self-intersection $-9$ and we obtain the configuration $C_7$ in $P_{12}$. Rationally blowing down $C_7$ gives a manifold $R_6$ homeomorphic to $P_6$.

\medskip

\hspace{.4in}\includegraphics[scale=0.50]{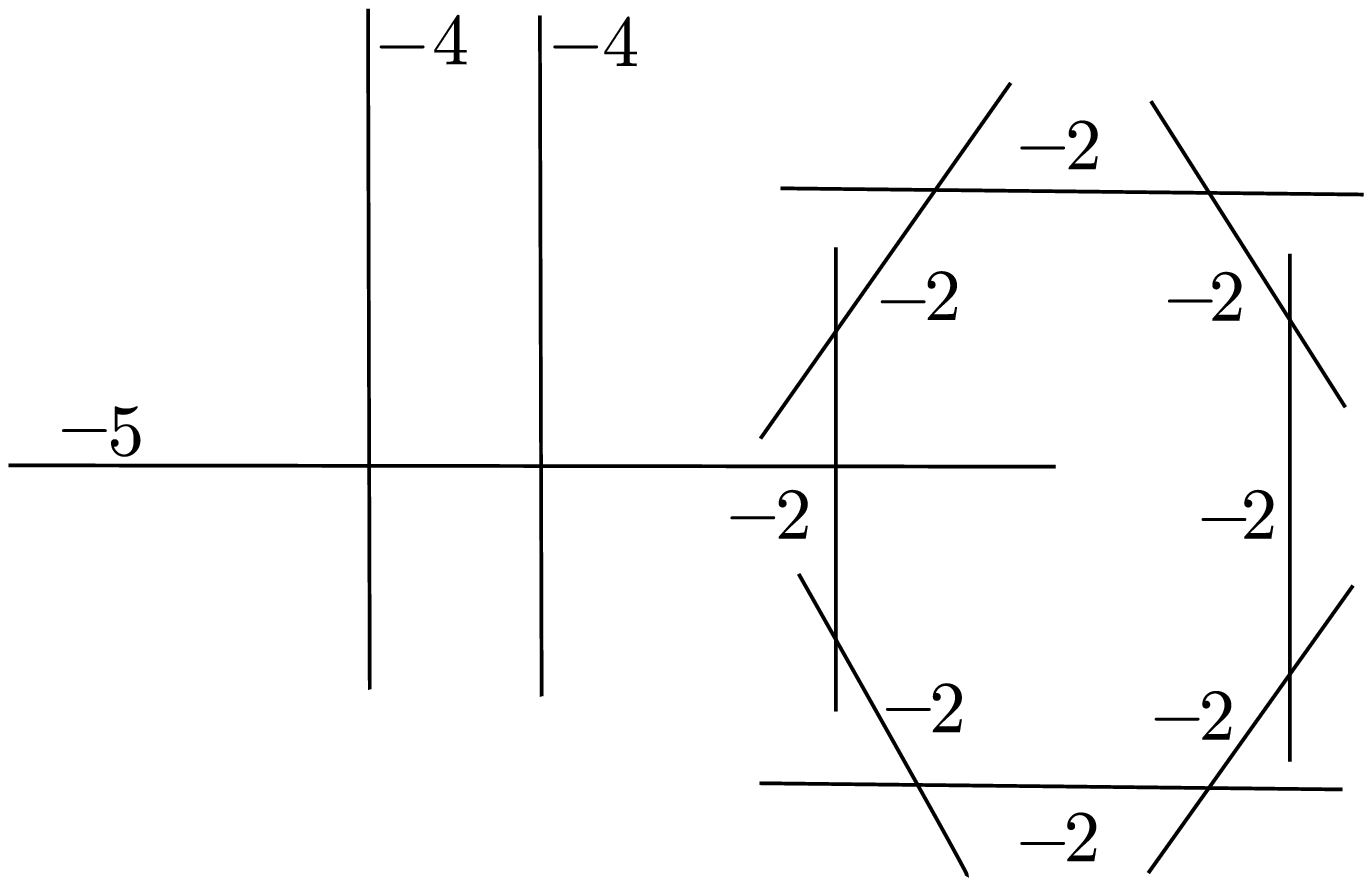}\\
\centerline{{Figure 5}}

\medskip

We shall work with the the first case, and then indicate what needs to be done to take care of the other cases. In Case 1, we obtain a manifold $R=R_8$ which is homeomorphic to $P_8$. We conjecture that $R$ is actually diffeomorphic to $P_8$, but for now it will suffice to see that it shares with $P_8$ the property that its Seiberg-Witten invariant vanishes. This will follow once we show that $R$ contains a sphere $H$ of self-intersection $1$ ({\it cf.} \cite{Turkey}). This sphere is obtained from the sphere $H$ representing a generator of $H_2(\CP;\Z)$ contained in $\CP\#9\CPb=P_9=E(1)$, because the construction of $R$ starts with $E(1)$, blows up and rationally blows down, and all the surfaces involved are disjoint from $H$. Hence:

\begin{prop} $\sw_R=0$.
\end{prop} 

Back in $E(1)=E(1)_K= E(1)\#_{F=S^1\x m} S^1\x M_K$ there is a nullhomologous torus $\L=S^1\x \lam$ where $\lam$ is the loop shown in Figure 6, which is a Kirby calculus depiction of $M_K = S^1\x S^2$, since $K$ is the unknot.

\noindent Since $\L$ as well as the $3$-manifold that it bounds ($S^1\x$ punctured torus) are disjoint from the regions where our constructions were made, $\L$ descends to a nullhomologous torus (which we still call $\L$) in $R$. Let $Q$ be the result of $0$-surgery on $\L \C R$, where the `$0$-framing' is taken from the $0$-framing on $\lam$ in Figure 6. After this surgery, the loop $\mu$ which bounds a normal disk to $\L$, does not bound in $Q$. In fact, $H_2(Q)$ is the direct sum of $H_2(R)$ with a hyperbolic pair generated by $\L_0$, the torus in $Q$ corresponding to $\L$, and a dual class represented by a torus built from the punctured torus that the longitude to $\lam$ (the surgery curve) bounds and the disk that the surgery curve bounds in $Q$. Thus $b^+(Q)=2$.

\vspace{.2in}
\hspace{1.4in}\includegraphics[scale=0.50]{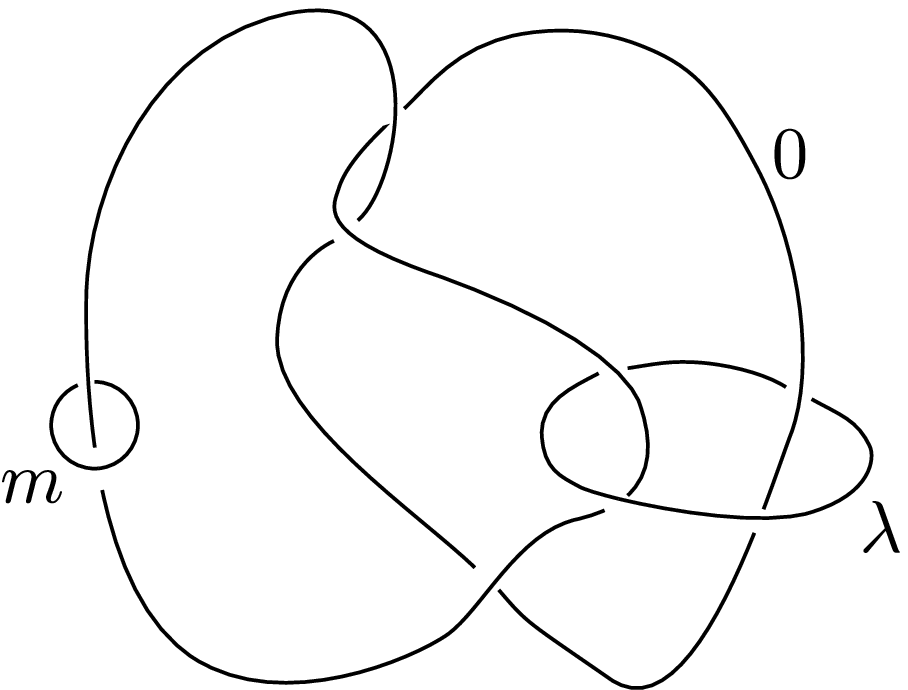}\vspace{.1in}\\
\centerline{{Figure 6}}

\smallskip

\begin{thm} The Seiberg-Witten invariant of $Q$ is: $\sw_{Q}=t^{-1}-t$, where $t=t_{S^1\x m}\in \Z H_2(Q)$.
\end{thm}
\begin{proof} The manifold $Q$ is obtained by:
\begin{itemize}\item[1.] Double node surgery with $K=$ the unknot, blowing up, then rationally blowing down.
\item[2.] $0$-surgery on $\L$.
\end{itemize}
Since $\L$ is disjoint from all the constructions in (1), the order in which (1) and (2) are performed is irrelevant. (As an aside, note that if we could exactly ``see'' $\L$ embedded in $P_8$, step (1) would be unnecessary, and we could then use $P_8$ rather than $R$.)

In $E(1)_K\cong E(1)$, do surgery on $\L$ first. Recall that $E(1)_K$ is a fiber sum 
\[ E(1)_K=E(1)\#_{F=S^1\x m} S^1\x M_K \]
and $M_K$ is the manifold given in Figure 6.  The result of $0$-surgery on $\L$ in $E(1)_K$ is the fiber sum
\[ E(1)_{K,0}=E(1)\#_{F=S^1\x m} S^1\x Y \]
where $Y$ is the $3$-manifold obtained from $0$-surgery on $\lam$ in Figure 6.

We shall now need the sewn-up link exterior construction of Brakes and Hoste. We recall what this means. Let
$L$ be a link in $S^3$ with two {\it{oriented}}\/ components $L_1$ and
$L_2$. Fix tubular neighborhoods $N_i\cong S^1 \x D^2$ of $L_i$ with
$S^1\x({\text{pt on $\bd D^2$}})$ a longitude of $L_i$, {\it i.e.}\/
nullhomologous in
$S^3\setminus L_i$. For any $A\in GL(2;\Z)$ with $\det A=-1$, we the get a
3-manifold
\[ s(L;A)=(S^3\setminus {\text{int}} (N_1\cup N_2))/ A \]
called a {\it{sewn-up link exterior}} by identifying $\bd N_1$ with $\bd N_2$ via a
diffeomorphism inducing $A$ in homology. For
$n\in\Z$, let
$A_n=\bigl(\begin{smallmatrix} -1&0\\-n&1 \end{smallmatrix}\bigr).$
A simple calculation shows that $H_1(s(L;A_n);\Z)=\Z\oplus\Z_{n-2\l}$ where
$\l\,$  is
the linking number in $S^3$ of the two components $L_1$, $L_2$, of $L$.
The second summand is generated by the meridian to either component.

We now refer to the proof of the Knot Surgery Theorem given in \cite{KL4M}. A key step in the proof, (see Figure 6 of \cite{KL4M}) which uses the work of Hoste \cite{H}, shows that $E(1)_{K,0}$ is diffeomorphic to 
\[ E(1)_L = E(1)\#_{F=S^1\x m} S^1\x s(L,A_{-2}) \]
where $L=L_1\cup L_2$ is the link of Figure 7, {\it{i.e.}} $L$ is the Hopf link. The orientations on $L_1$ and $L_2$ are inherited from fixing an orientation on the knot $K$, {\it{e.g.}} in Figure 6. Note that since the linking number of $L_1$ and $L_2$ is $-1$, we have
$H_1(s(L;A_{-2});\Z)=\Z\oplus\Z$.  

\smallskip

\hspace{1.4in}\includegraphics[scale=0.50]{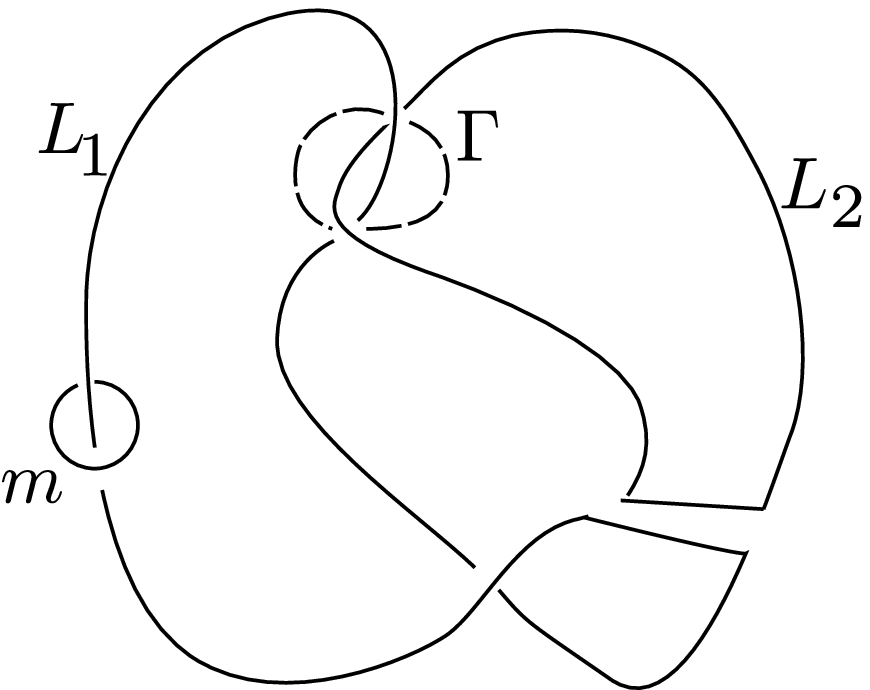}\\
\centerline{{Figure 7}}

\smallskip

In $s(L,A_{-2})$ we see an embedded torus transverse to the meridian $m$ obtained by sewing up a Seifert surface for $L$, and we also see a loop $\G$ which satisfies conditions ({\em i}) and  ({\em ii}) for double node surgery.

The proof of the Knot Surgery Theorem shows that the Seiberg-Witten invariant of $E(1)_L$ can be calculated via skein moves (macarena). Note that we are now calculating the Seiberg-Witten invariant of a manifold with $b^+=2$. This calculation is shown in Figure 8. This figure depicts the fact that 
\[\sw_{E(1)_L} =\sw_{E(1)_U} - (t-t^{-1})^2\sw^-_{E(1)_{K_0}} \]
where $U$ is the unlink and $K_0$ is the unknot. (See Section 3, equation (3) of \cite{KL4M}.)

A $2$-sphere that separates the two components of the unlink, gives rise to an essential $2$-sphere in $E(1)_U$ of self-intersection $0$. Thus $\sw_{E(1)_U} =0$. Since $E(1)_{K_0}=E(1)$, we have $\sw^-_{E(1)_{K_0}} =(t-t^{-1})^{-1}$. Thus
$\sw_{E(1)_L}=t^{-1}-t$. This completes our discussion of step (2).

\vspace{.2in}
\hspace{.3in}\includegraphics[scale=0.45]{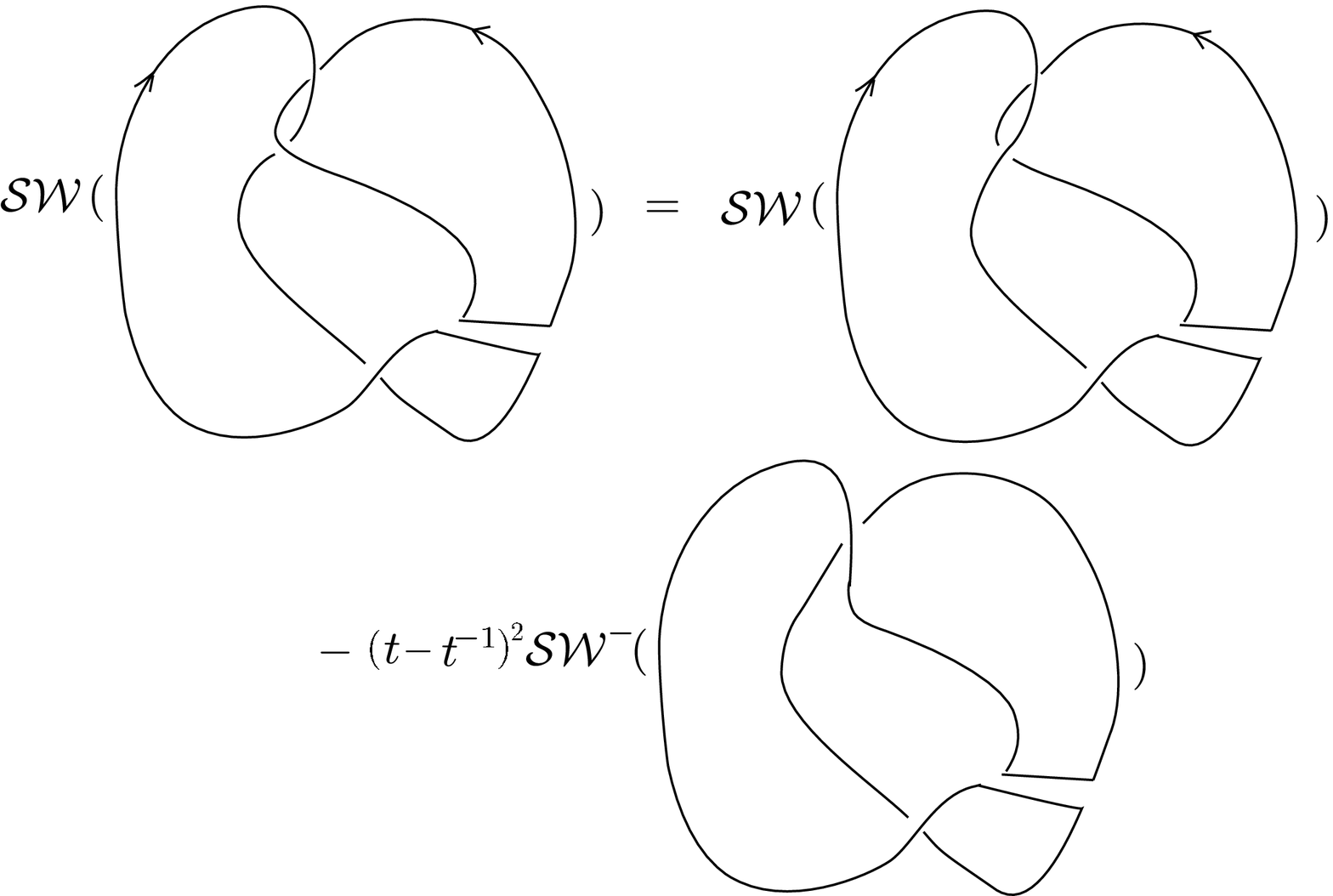}\vspace{.1in}\\
\centerline{{Figure 8}}

Next we carry out the constructions of step (1). In $E(1)_L$ there is a genus one pseudosection to which we can apply the double node construction. This pseudosection is the connected sum of a section in $E(1)$ with the torus in \{point\}$\x s(L,A_{-2})$ obtained by sewing up the shaded region in Figure 9. The necessary loop $\G$ is shown in Figure 7.

\vspace{.2in}
\hspace{1.4in}\includegraphics[scale=0.50]{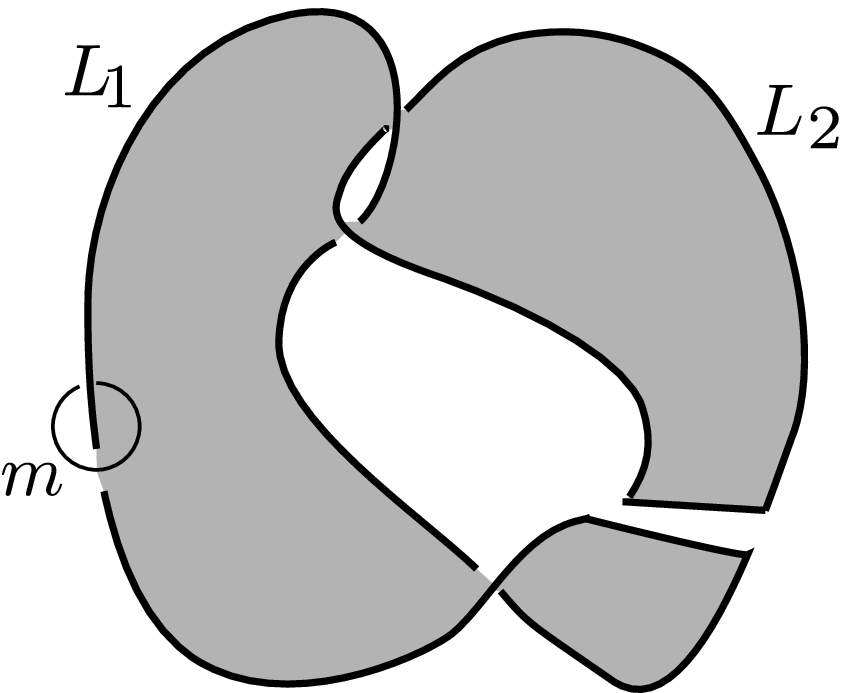}\\
\centerline{{Figure 9}}

The result of the double node construction is an immersed genus $0$ pseudosection with one positive double point. Blow up at this double point to get an embedded $2$-sphere $C$ of self-intersection $-5$ in $E(1)_L\#\CPb$. The blowup formula \cite{Turkey} implies that $\sw_{E(1)_L\#\CPb}=(t^{-1}-t)(e+e^{-1})$ where $e$ is is the class in the group ring corresponding to the exceptional curve. Hence the basic classes of $E(1)_L\#\CPb$ are $\pm F\pm E$. Now $\pm(F+E)\cdot C=\pm 3$ whereas 
$\pm(F-E)\cdot C=\mp 1$. It follows from the Rational Blowdown Theorem that only the basic classes $\pm(F+E)$ descend to the rational blowdown, $Q$. Thus $Q$ has two basic classes whose Seiberg-Witten invariants are those of $\pm(F+E)$ in $E(1)_L\#\CPb$, namely, $\mp1$. 
\end{proof}

\begin{thm} There are infinitely many nonpairwise diffeomorphic $4$-manifolds homeomorphic to $P_8=\CP\#8\CPb$ obtained from $1/n$-surgery on the nullhomologous torus $\L$ in $R$.
\end{thm}
\begin{proof} For $n\ge 2$ let $X_n$ be the $4$-manifold obtained from $1/n$-surgery on $\L$ in $R$. By this we mean $S^1$ times $1/n$-surgery on $\lam$ in Figure 6. It is easy to see that $X_n$ is simply connected and that there is an isomorphism $\vp: H_2(X_n;\Z)\to H_2(R;\Z)$, which is realized outside of a neighborhood of the surgery by the identity map. Thus $X_n$ is homeomorphic to 
$P_8$. 

Morgan, Mrowka, and Szabo have calculated the effect of such a surgery on Seiberg-Witten invariants \cite{MMS}.  Given a class $k\in H_2(X_n)$:
\[ \ssw_{X_n}(k)=\ssw_R(\vp(k')) + n\,\sum_i\ssw_{Q}(k''+i[\L_0])\]
(Recall that the torus $\L$ is nullhomologous in $R$ and the corresponding torus $\L_n$, the core of the surgery, is nullhomologous in $X_n$.)  Further, $k''\in H_2(Q)$ is any class which agrees with the restriction of $k$ in $H_2(R\- \L\x D^2,\bd)$ in the diagram:
\[ \begin{array}{ccc}
H_2(Q) &\longrightarrow & H_2(Q, \L_0\x D^2)\\
&&\Big\downarrow \cong\\
&&H_2(R\- \L\x D^2,\bd)\\
&&\Big\uparrow \cong\\
H_2(X_n)&\longrightarrow & H_2(X_n,\L_n\x D^2)
\end{array}\]
The Seiberg-Witten invariants of the two $b^+=1$ manifolds $X_n$ and $R$ are calculated in corresponding chambers. 

Given $k\in H_2(X_n)$ and $H$ an element of positive self-intersection in $H_2(X_n)$, the small perturbation chamber, {\it i.e.} the sign $\pm$ such that $\ssw_{X_n}(k)=\ssw_{X_n,H}^{\pm}(k)$ is determined homologically. This means that the small perturbation chambers for $k$ in $X_n$ and for $\vp(k)$ in $R$ correspond under $\vp$. According to the previous theorem, there are only two classes, $\pm T$, $T=[S^1\x m]$ in $Q$ with nontrivial Seiberg-Witten invariants, and $\sw_{Q}(\pm T)=\mp 1$.

Thus we have 
\[ \sw_{X_n} = \sw_R +\sw_{Q} = 0 + n\,(t^{-1}-t) \]
This shows that the manifolds $X_n$ are pairwise nondiffeomorphic. That they are all minimal follows from the blowup formula.
\end{proof}

\section{Infinite families with $b^-=5,6,7$}

For $b^-=6,7$, the constructions and calculations that there is a nullhomologous torus $\L$ in $R=R_6$
or $R_7$ such that $1/n$-surgery on $\L$ in $R$ produces infinitely many distinct smooth structures on $P_6$
or $P_7$ 
 are completely analogous to those in the last section. What differs is the proof that the Seiberg-Witten invariants of $R=R_6$
or $R_7$ vanish. We shall accomplish this by using an argument adapted from \cite{OS}. This has in common with the previous argument its dependence on the adjunction inequality.
The important point in the argument is that we are starting our construction with $E(1)_K=E(1)$; so all the exceptional curves are represented by spheres, etc. First consider the $b^-=7$ case. The classes
\begin{eqnarray*} V_1  = H-E_1-E_2-E_3, &\, \, &V_2  =  H-E_2-E_3-E_4\\ V_3  =  H-E_3-E_4-E_5,  && V_4 = H-E_6-E_7-E_8-E_9\\ V_5 = F-E_5,  && V_6 =E_{11}-E_1-E_2\\  V_7 = E_{10}-E_1-E_2, && V_8= 2H-3E_{11}\end{eqnarray*}
are all orthogonal to the configuration $C_5$ and generate $H_2(P_{11}\- C_5;\Z)=H_2(R_7\- B_5;\Z)$. In $P_{11}$, the classes $V_1, V_2, V_3$ are represented by embedded spheres of self-intersection $-2$, $V_4, V_6, V_7$ are represented by embedded spheres of self-intersection $-3$, $V_5$ is represented by an embedded torus of self-intersection $-3$, and $V_8$ is represented by an embedded torus with square $-5$. According to the argument of \cite{OS}, any basic class $k$ of $R_7$ must satisfy the adjunction inequality 
\[ V_i^2 + |k\cdot V_i| \le 0, \ i=1,\dots,8 \]
Furthermore, $k$ must satisfy, 
\[ k^2\ge 2, \hspace{.25in} k^2\equiv 2\pmod8  \]
These follow from the fact that for any basic class, its corresponding moduli space must  have nonnegative even dimension. 

Since $H^2(R_7,\Z)$ injects into $H^2(R_7\- B_5;\Z)$, any Seiberg-Witten basic class of $R_7$ is uniquely determined by its intersection numbers with $V_1,\dots , V_8$. There is now a finite check for possible basic classes $k$ of $R_7$ which must satisfy these three previous conditions. This check turns up $40$ classes in $H^2(R_7\- B_5;\Z)$. Another class in $H_2(P_{11};\Z)$ which is orthogonal to $C_5$ is $V_9=2H-3E_{10}$. It is represented by an embedded torus of self-intersection $-5$. Any basic class of $R_7$ must also satisfy the adjunction inequality with respect to $V_9$. This condition reduces the number of possible classes to $14$. According to the Rational Blowdown Theorem the Seiberg-Witten invariant of any such class is determined by the Seiberg-Witten invariant of an appropriate lift to $P_{11}$. Such a lift determines a Seiberg-Witten moduli space for $P_{11}$ whose formal dimension is the same as that of the moduli space for $R_7$ corresponding to the class being lifted. This is accomplished via an extension across $C_5$ for each of the $14$ possibilities. 

The class $\hat{H}=V_1+V_2+V_3+V_8$ is orthogonal to $C_5$ and $\hat{H}^2=4>0$, $\hat{H}\cdot H=6>0$. Hence $\hat{H}$ serves as a period point for $R_7$. Since $\sw_{P_{11},H}=0$, for any characteristic cohomology class $k$ of $P_{11}$, 
\[ \sw_{P_{11},\hat{H}}(k)=\begin{cases} 0 \ \ &{\text{if the signs of $k\cdot \hat{H}$ and $k\cdot H$ agree}}\\
\pm 1 & {\text{if the signs of $k\cdot \hat{H}$ and $k\cdot H$ do not agree}}
\end{cases}  \]
Using this criterion on each of the $14$ possibilities mentioned above shows that each has Seiberg-Witten invariant equal to $0$.
Thus we have $\sw_{R_7}=0$. 

The $b^-=6$ case follows similarly, but the calculation turns out to be easier.  The classes 
\begin{multline*} V_1=E_{12}-E_1,\, V_2=E_{11}-E_1\, V_3=E_{10}-E_1, \, V_4=E_3-E_1, \\ V_5=E_2-E_1, \, V_6=H-3E_1, \, V_7=2F+H-E_3-E_{10}-E_{11}-E_{12} \end{multline*}
generate $H_2(P_{12}\- C_7;\Z)=H_2(R_6\- B_7;\Z)$. In $P_{12}$, the classes $V_1, \dots, V_5$ are all represented by embedded spheres of self-intersection $-2$, $V_6$ is represented by an embedded torus of self-intersection $-8$, and $V_7$ is represented by an embedded surface of genus $4$ with square $5$.  If $k$ is a basic class of $X$, then b it must satisfy the adjunction inequality with respect to the classes $V_1,\dots,V_7$; {\it{i.e.}} 
\[ V_i^2 + |k\cdot V_i| \le 0, \ i=1,\dots,6;  \  \ V_7^2 + |k\cdot V_i| \le 2g-2=6\]
and as above, $k$ must also satisfy
\[  k^2\ge 3, \hspace{.25in} k^2\equiv 3\pmod8  \]
This time a check turns up no classes in $H^2(R_6\- B_7;\Z)$ which satisfy these conditions.
Hence $\sw_{R_6}=0$.

To obtain families of manifolds homeomorphic to $P_5$, we start with an elliptic fibration on $E(1)$ with two nodal fibers, two $I_2$-fibers and an $I_6$-fiber. (Again, see \cite{Pn}.)
This time we have two double node neighborhoods. One goes through the same construction in each of these to obtain an immersed pseudosection with two double points. Blowing up each, we obtain a sphere of square $-9$ in $E(1)\# 2\CPb$, and adding on five of the spheres of the $I_6$-singularity, gives a copy of the configuration $C_7$, which can be rationally blown down to obtain a manifold $R_5$ homeomorphic to $P_5$. This process can be carried out so the the $+1$-sphere $H$ descends to $R_5$; so we see that $\sw_{R_5}=0$. Furthermore, we get the nullhomologous tori $\L_1, \L_2$, as before; one in each neighborhood. Performing $0$-surgery on each gives a manifold $Q$ with $b^+=3$ and $\sw = (t_1^{-1}-t_1) (t_2^{-1}-t_2)$. Perform $+1$-surgery on $\L_2$ and $1/n$-surgery on $\L_1$ to obtain $Y_n$ which is homeomorphic to $P_5$, but which has $\sw_{Y_n}=n (t_1^{-1}-t_1) (t_2^{-1}-t_2)$.

\end{document}